\documentclass[11pt, reqno]{amsart}
\usepackage{geometry}                
\usepackage{graphicx}
\usepackage{amssymb}
\usepackage{epstopdf}
\usepackage{color}

\usepackage{caption}
\usepackage{subcaption}
\usepackage{pdfpages}

\usepackage{mathrsfs}
\usepackage{enumitem}
\usepackage{hyperref}
\usepackage{isomath}

\numberwithin{equation}{section}

\title[Dilogarithm identities ]
{Dilogarithm identities after Bridgeman
}

\keywords{Rogers dilogarithm, identities, Fibonacci numbers, Lucas numbers, continued fractions}
\author{Pradthana Jaipong, Mong Lung Lang,  Ser Peow Tan, and Ming Hong Tee}
\date{}                                           

\address{Department of Mathematics, Chiang Mai University\\
	Thailand } \email{pradthanadee@gmail.com}
\address{
	Singapore } \email{lang2to46@gmail.com}
\address{Department of Mathematics \\
	National University of Singapore \\
	Singapore 119076} \email{mattansp@nus.edu.sg}
\address{Department of Mathematics \\
	National University of Singapore \\
	Singapore 119076} \email{e0174864@u.nus.edu}

\subjclass[2000]{}

\thanks{The third author is partially supported by the National University
	of Singapore academic research grant R-146-000-289-114. The first and third author are grateful to the Temasek foundation for support.}   

\date{\today}

\begin{document}
\newtheorem{theorem}{{\bf Theorem}}[section]
\newtheorem{lemma}[theorem]{{\bf Lemma}}
\newtheorem{proposition}[theorem]{{\bf Proposition}}
\newtheorem{corollary}[theorem]{{\bf Corollary}}
\newtheorem{definition}[theorem]{{\bf Definition}}
\newtheorem{example}[theorem]{{\bf Example}}
\newtheorem{remark}[theorem]{{\bf Remark}}

\newcommand{\LL}{\mathcal{L}}
\newcommand{\ZZ}{\mathbb{Z}}
\newcommand{\RR}{\mathbb{R}}
\newcommand{\CC}{\mathbb{C}}
\newcommand{\NN}{\mathbb{N}}
\newcommand{\HH}{\mathbb{H}}

\vspace{-0.3in}

\begin{abstract}
	Following Bridgeman, we demonstrate several families of infinite dilogarithm identities associated with Fibonacci numbers, Lucas numbers, convergents of continued fractions of even periods, and terms arising from various recurrence relations.

\end{abstract}

\maketitle
\tableofcontents

\section{Introduction}

In this article, we exhibit several families of infinite identities involving the Rogers dilogarithm $\LL(x)$, following Bridgeman. These identities arise from Bridgeman's orthospectral identity (see \cite{B1}) as applied to various hyperbolic cylinders. They generalize the connection found by Bridgeman in \cite{B2} between the solutions of Pell's equations, the continued fraction convergents of these solutions and the Rogers dilogarithm. In particular, families of identities involving the dilogarithms of Fibonacci numbers, Lucas numbers, other recurrence sequences and convergents of continued fraction expansions with  period
 two or even period are derived.
\subsection{Main results}
 Recall that the Rogers dilogarithm $\LL(z)$, for $0 \le z \le 1$ is given by
\[\LL(z)=Li_2(z)+\frac{1}{2}\log|z|\log (1-z), \quad {\hbox{where}} \quad Li_2(z)=\sum_{n=1}^\infty \frac{z^n}{n^2}, \qquad |z|\le 1,\] is the dilogarithm function. Note that $\LL(0)=0$ and $\LL(1)=\pi^2/6$ and although $\LL(z)$ can be extended by analytic continuation to complex values, we will be mainly concerned with its values for $z \in [0,1]$. The dilogarithm and the Rogers dilogarithm appear in various forms in Algebraic K-theory, mathematical physics, number theory and hyperbolic geometry, see for example \cite{Zag}. The Fibonacci numbers are defined by $f_0 =1, f_1=1,$  $f_n=f_{n-1}+f_{n-2}$, $n \in \ZZ$ and the Lucas numbers are defined by
$l_0 = 2, l_1=1, l_n= l_{n-1}+l_{n-2}$, $n \in \ZZ$ and $\phi := (1+\sqrt5)/2$ is the golden ratio. 

\medskip

\noindent {\it Remark}. In the statement of the results in the rest of the introduction, the equation number (X.Y) attached to the identities indicate the section where it is stated and proved, for example (8.7) indicates that the identity is proven in section 8 as equation (8.7).

\begin{theorem} Let $\LL(x)$ be the Rogers dilogarithm, $ f_k$, $l_k$ the Fibonacci and Lucas numbers and $\phi$ the golden ratio. Let $\epsilon_1=1/2$ and $\epsilon_k=1$ otherwise. We have:

$$\sum_{k=2}^{\infty }
\LL \left (\left (
\frac{f_{2n}}{f_{2nk}}\right )^2\right )=\LL \left( \left (\frac{1}{\phi}\right )^{4n}\right ), \qquad n \in \mathbb{N}. \eqno(8.7)$$

When $n=1$,

$$ \sum_{k=2}^{\infty}
\LL  \left (\left ( \frac{ 1}
{f_{2k}}\right )^2\right )=\LL (1/\phi^4 ) ,\qquad
\sum_{k=1}^{\infty } 
\LL \left (\frac{1}{f_{2k-3}f_{2k+1}}\right ) =
\LL (1-1/\phi^4),\eqno(5.8)$$
$$\sum_{k=1}^{\infty }
\left ( \LL \left (\frac{1}{f_{2k+2}^2}\right )
+ \LL \left (\frac{1}{f_{2k-3}f_{2k+1}}\right )\right ) = \frac{\pi^2}{6}.\eqno(4.7)$$
When $n=2$,

$$
\sum_{k=2}^{\infty }
\LL \left (\frac{3^2}{f_{4k}^2}\right )=\LL (1/\phi^8)
,\qquad
\sum_{k=0}^{\infty }
\LL \left (\frac{45}{l_{4k-2} l_{4k+6}}\right)=\LL (1-1/\phi^8),
\eqno(5.9)$$

 $$\sum_{k=2}^{\infty }
 \LL \left (\frac{3^2}{f_{4k}^2}\right)
 +
 \sum_{k=0}^{\infty }
 \LL \left (\frac{45}{l_{4k-2}l_{4k+6}}\right)
 = \frac{\pi^2}{6}.\eqno(4.8)$$
 
 \medskip
 
 For powers of $\phi$ congruent to $2 \mod 4$, we have:

  $$\sum_{k=2}^{\infty} \LL  \left ( \frac{l_{2n+1}^2}{l_{k(4n+2)-(2n+1)}^2}\right )
 +
 \sum_{k=1}^{\infty} \LL  \left ( \frac{l_{2n+1}^2}{5f_{k(4n+2)}^2}\right )=\LL  (1/\phi^{4n+2}), \qquad n \in \mathbb{N}\cup\{0\}.
 \eqno(10.1) $$

\medskip

We also get identities for $\pi^2/6$ and $\pi^2/10=\LL(1/\phi)$ where the arguments of the terms in the infinite sums are expressed in terms of the Fibonacci numbers and the Lucas numbers:

$$ \sum_{k=1}^{\infty } \left (
\LL \left (\frac{1}{5f_{2k}^2}\right)
+
\LL \left (\frac{1}{l_{2k+1}^2}\right)
+
\LL \left (\frac{1}{l_{2k-2} l_{2k}}\right)
+
\epsilon_k   \LL \left (\frac{1}{5f_{2k-3}f_{2k-1}}\right)\right )
= \pi^2/6
,\eqno(12.2)$$


$$\sum_{k=1}^{\infty }\left (
\LL \left (\frac{1}{l_{2k-2} l_{2k}}\right)
+
\epsilon_k \LL \left (\frac{1}{5f_{2k-3}f_{2k-1}}\right)\right )
=\LL (1/\phi)= \frac{\pi^2}{10}.\eqno(12.3)$$

 \end{theorem}

An identity for $\pi^2/12=\LL(1/2)$ is given below in Corollary 1.4.  Identities  for $\LL (1/\phi^{2n+1})$ are more involved and can be found in Section \ref{s:twokplusone}.
 Note that (10.1) for $n=0$   was first derived by Bridgeman in \cite{B2}.
 
 \medskip
 
The next result is a generalization of the Richmond Szekeres identity \cite{RS}, see also \cite{Le}. Define the cross ratio of $4$  points in $\hat{\mathbb{C}}$, at least three of which are distinct, by
\[ [z_1,z_2,z_3,z_4]=\frac{(z_1-z_2)(z_4-z_3)}{(z_1-z_3)(z_4-z_2)}.\]
\begin{theorem}
Consider the ideal hyperbolic polygon $P$ with  vertices $v_1,\ldots, v_{n}$, $n \ge 3$ with  $ v_1=0/1 < v_2 <\ldots< v_{n-1}=1/1$, $v_n=\infty$.
Then

\[
\mathop{\sum\sum}_{1 \le j<i\le n-2}
\LL \left(
[v_i,v_{i+1},v_{j},v_{j+1}]
\right)+ \sum_{i=1}^{n-2}\sum_{j=1}^{n-2} \sum_{k=1}^{\infty}
\LL \left(
[v_i,v_{i+1},k+v_{j},k+v_{j+1}]
\right)
=\frac{(n-2)\pi^2}{3}.\\  \eqno(13.2)   
\]

\end{theorem}

\noindent {\it Remark:} The case $n=3$ reduces to the Richmond Szekeres identity  $\sum_{k=2}^{\infty} \LL  \left (\dfrac{1}{k^2}\right ) = \pi^2/6$. This case also follows as limits of (5.3) or (7.1) below. More interesting examples with $n=4$ can be found in Section 13. Furthermore, note that if adjacent vertices of $P$ are Farey neighbors, then the arguments of $\LL(x)$ in the identity are all rational numbers with numerator $1$.

\medskip

We also have the following identities involving sequences defined by recurrences in the next two results:

\begin{theorem} (Theorem 5.1)
	 Suppose that $t>2$ and $u>1/u$ are the roots of $x^2-tx+1=0$. Then $t=u+1/u$ and
		$$
		\sum_{n=1}^{\infty }
		\LL \left ( \frac {1}{q_{n}\,^2}\right )
		=\LL (1/u^2),\qquad
		\sum_{n=1}^{\infty }
		\LL  \left (\frac{t-2}{(q_{n}-q_{n-1})(q_{n-2}-q_{n-3})}
		\right)=\LL ( 1-1/u^2),
		\eqno(5.3)$$
		
\noindent where $\{q_n\}$ is the recurrence defined by
		$  q_0=1 , q_1=t,
		q_n= tq_{n-1}-q_{n-2}.$

\end{theorem}

Note that $\LL(1/u^2) +\LL(1-1/u^2)= \pi^2/6$ so the two identities above can be combined to give an expression for $\pi^2/6$.
In the case $t = \sqrt n+1/\sqrt n$, Theorem 1.3 gives the following
 identity.

\begin{corollary}
Let $n>2$ be an integer. Then
$$\mathcal
L\left (\frac{1}{n}\right )
= \sum_{k=1}^{\infty} \LL \left (\left (
\frac {n^{k/2}}{ n^k +n^{k-1}+\cdots +n^2+n+1}\right ) ^2\right ).\eqno(7.1)$$
In particular, when $n=2$, we get
$$\frac{\pi^2}{12}=\mathcal
L\left (\frac{1}{2}\right )
= \sum_{k=1}^{\infty} \LL \left (
\frac {2^k}{ (2^{k+1}-1)^2} \right ).$$

\end{corollary}

\begin{theorem}
	Let
	$A={\tiny
		\left ( \begin{array}{cc}
		a& c\\
		b & d\\
		\end{array}
		\right ) }\in SL(2, \Bbb Z)$ with positive entries, trace $t=a+d>2$ and eigenvalues $u>1/u$. Set
	$A^{n}=
	{\tiny \left ( \begin{array}{cc}
		p_{2n-1}& p_{2n-2}\\
		q_{2n-1} & q_{2n-2}\\
		\end{array}
		\right )} $, where
	$A^{0}=
	{\tiny \left ( \begin{array}{cc}
		p_{-1}& p_{-2}\\
		q_{-1} & q_{-2}\\
		\end{array}
		\right ) =
		\left ( \begin{array}{cc}
		1& 0\\
		0 & 1\\
		\end{array}
		\right )} .$
	Then	 $A^n= tA^{n-1}-A^{n-2}$ and
 $$
         \sum _{n=2}^{\infty}
\LL \left (\left (\frac {b}{q_{2n-1}}\right)^2 \right )=\LL (1/u^2),
\eqno(8.5)$$

$$
2\sum _{n=2}^{\infty}
\LL \left (\left (\frac {b}{q_{2n-1}}\right)^2 \right )
+\sum_{n=1}^{\infty }\LL \left (
\frac{bc}{q_{2n}q_{2n-4}}\right ) +
\sum_{n=1}^{\infty }\LL \left (
\frac{bc}{p_{2n+1}p_{2n-3}}\right )
+\LL \left (\frac{bc}{ad}\right ) = \pi^2/3.\eqno(8.6)
$$

\end{theorem}

\medskip
Theorem 1.5 can be applied to periodic continued fractions of period two to obtain the following corollary.

\begin{corollary}
	Let $\alpha=[\overline {a,b}]$ be a continued fraction of
	period $2$ and let $r_n = p_n/q_n$ be the $n$-th convergent of $\alpha$.
	Set
	$A=  {\tiny
		\left ( \begin{array}{cc}
		ab+1&  a\\
		b &1\\
		\end{array}
		\right )}
	=
	{\tiny
		\left ( \begin{array}{cc}
		p_1&  p_0\\
		q_1&q_0\\
		\end{array}
		\right )
	},$
	with eigenvalues $u>1/u$. Then
	$$ \LL (1/(b\alpha+1)^2) =
	\LL (1/u^2) =          \sum _{n=2}^{\infty}
	\LL \left (\left (\frac {b}{q_{2n-1}}\right)^2 \right ).
	\eqno(9.2)$$

\end{corollary}

The case $b=1$ above was first derived by Bridgeman in \cite{B2}. For the more general case of  periodic continued fractions of even period greater than two, the identity is more involved and is studied in Section 15. The precise statement for the identities is stated as Theorem 15.3. The crowns in these cases are more interesting as the number of tines will be greater than one.

%

\medskip

\subsection{Geometric background}
\begin{figure}[h]
\begin{center}
	\scalebox{0.15}[0.15]{\includegraphics{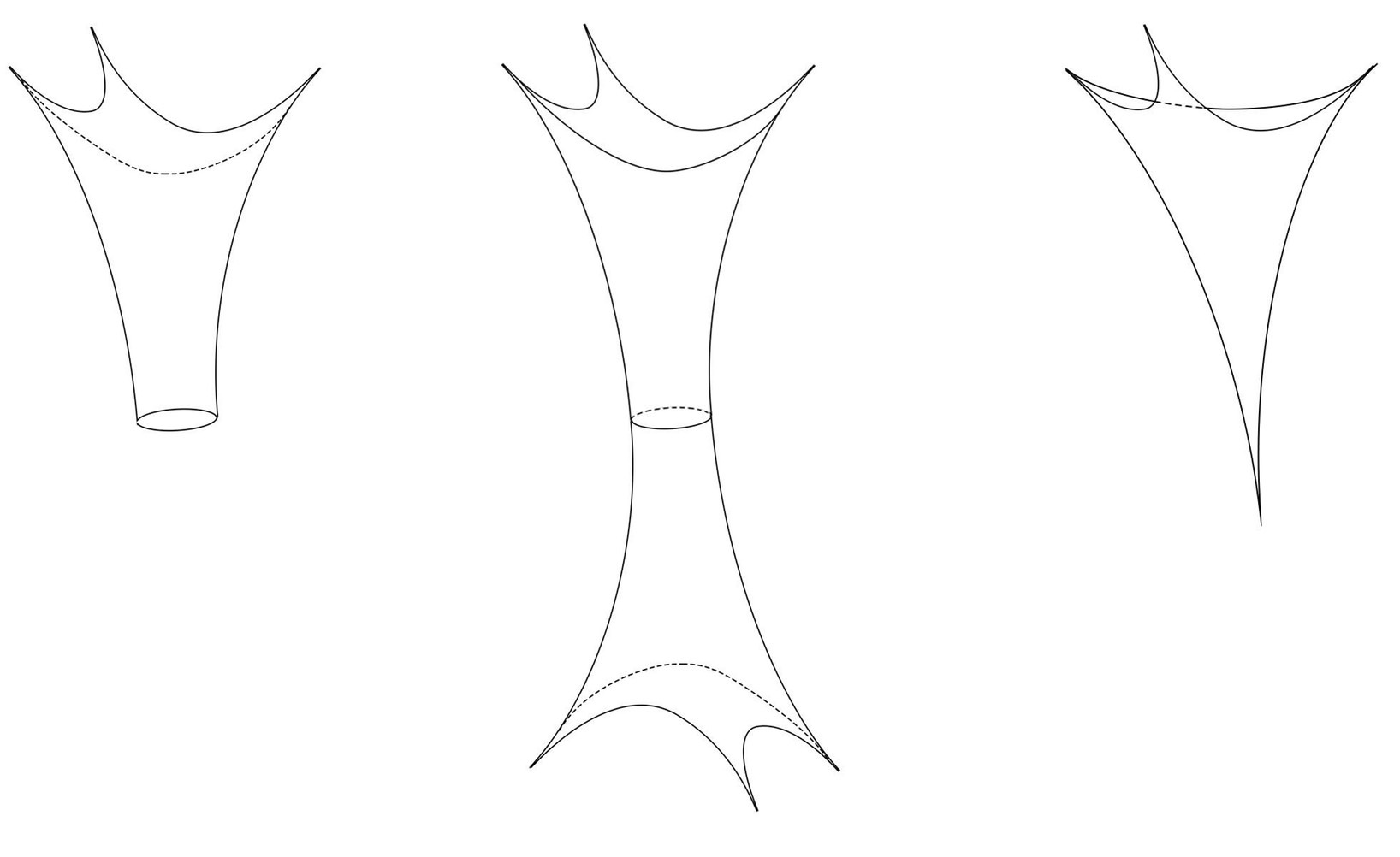}}
\end{center}

\caption{A crown; double crown; and degenerate crown with parabolic holonomy} \label{fig:crowns}
\end{figure}

The geometry behind our identities is simple, the geometric objects we consider are just hyperbolic cylinders $S$ with finite area, which have infinite cyclic fundamental group $\langle g\rangle$, with holonomy $\rho(g):=T \in  PSL (2, \RR)$. These cylinders  have a finite number of  boundary components which are either closed geodesics or complete infinite geodesics adjacent to {\it boundary cusps} (as opposed to regular cusps), and are called crowns or double crowns, see figure \ref{fig:crowns}. The boundary cusps are the tines of the (double) crowns. Typically, $T$ is hyperbolic, the degenerate situation where $T$ is parabolic (in which case there is one regular cusp) gives rise to Theorem 1.2. The geometry will be encoded in what we call {\em feasible pairs} $(T,P)$ where $P$ is a hyperbolic polygon with vertices in $\HH \cup \partial \HH$ and $T=\rho(g)\in PSL (2, \Bbb R)$ identifies two sides of $P$ and maps $P$ to the exterior of $P$, so that $P$ is a fundamental domain for the cylinder $S$. Judicious choices of the pair $(T,P)$ relating the geometry to the algebra and arithmetic of various sequences and application of Bridgeman's orthospectral identity  gives us the identities stated above.


\begin{figure}[h]
	\begin{center}
		\scalebox{0.3}[0.3]{\includegraphics{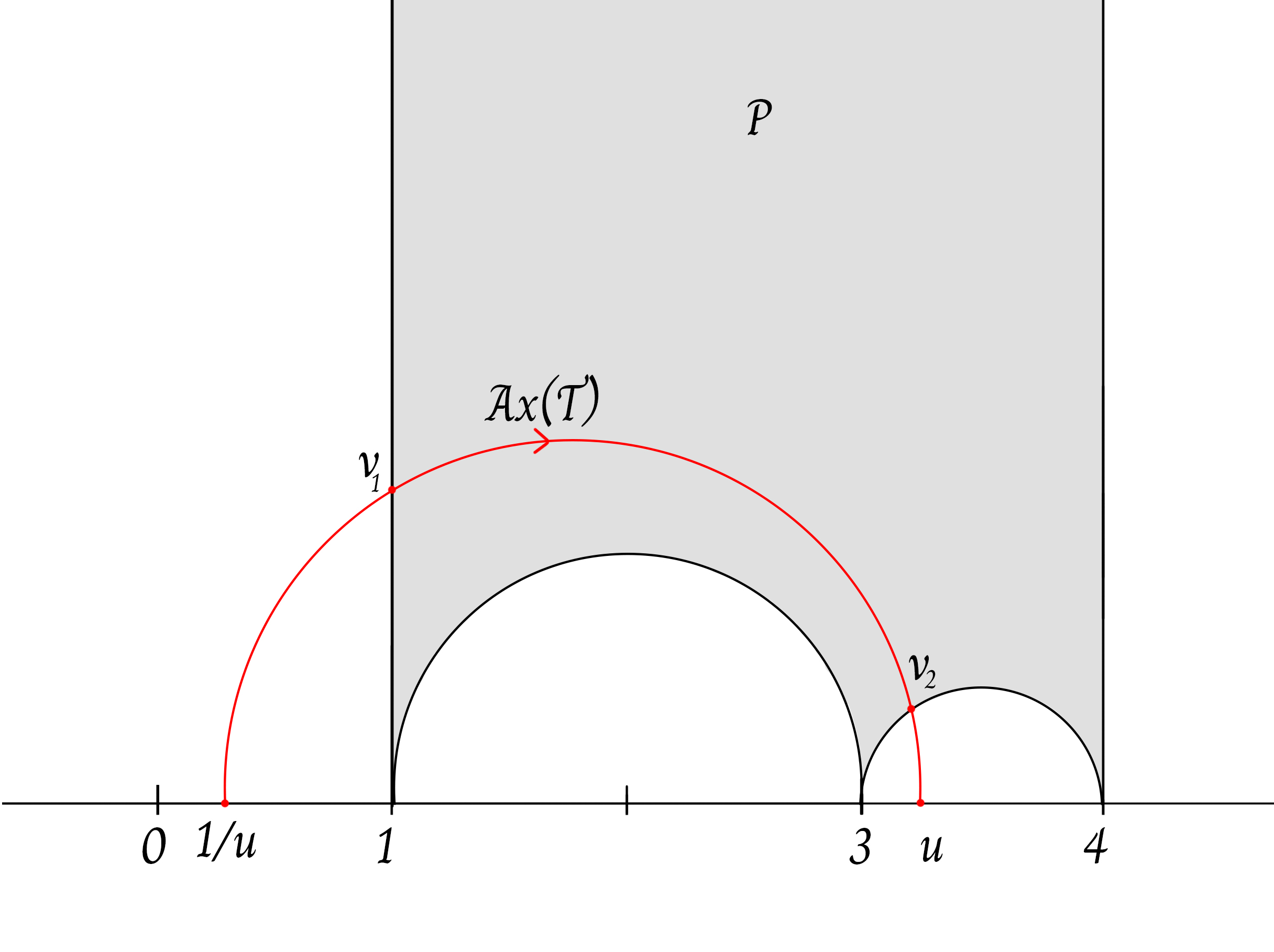}}
	\end{center}
	
	\caption{A feasible pair showing $P$ and the axis of $T$, fixed points $1/u$ and $u$ of $T$, and intersection points $v_1,v_2$.} \label{fig:FP}
\end{figure}


\subsection{The feasible pairs}\label{ss:feasible}
We call a pair $(T, P)$ {\em feasible} if $T\in PSL (2, \Bbb R)$
identifies two  sides of a hyperbolic polygon $P$ (the hyperbolic convex hull of a finite number of vertices in $\HH \cup \partial \HH$)
and sends the interior of $P$ to the exterior of $P$. In the case where $P$ has a finite side, than this side must be part of the invariant axis for $T$. For convenience, we will often choose a lift of $T$ to $SL (2, \Bbb R)$ with non-negative trace which we will also denote by $T$ and work with this lift.
All our  identities are constructed by studying
such feasible pairs (models).  Let $t$ be the trace of $T$.
The (degenerate) model $(T, P)$ where $t=2$ is studied in Section 13.
In the case $t>2$, we study the following  three feasible pairs giving rise to double-crowns. A modification (see Section 5) gives related pairs $(T,P')$ giving rise to crowns. The polygon $P$ for continued fractions of even period $>2$ is more complicated, with more tines, see Section 15.
\begin{enumerate}
	\item[(i)]
	$T ={\tiny\left ( \begin{array}{cc}
		t&  -1\\
		1 &0\\
		\end{array}
		\right )}$ and $P$ is the hyperbolic convex hull of $\{ 1, t-1, t ,\infty\}$ where $t >2$, see figure \ref{fig:FP} for the case $t=4$.

	\item[(ii)]
	$T ={\tiny\left ( \begin{array}{cc}
		a&  c\\
		b &d\\
		\end{array}
		\right )}$, $a, b , c, d >0$  and $P$ is the hyperbolic convex hull of $\{ c/d, a/b, t ,\infty\}$ where $a+d>2$.
 This feasible pair works better then (i)  for continued fractions.
	
	\item
	[(iii)]
	$T ={\tiny\left ( \begin{array}{cc}
		t&  -1\\
		1 &0\\
		\end{array}
		\right )}$ and $P$ is the hyperbolic convex hull of $\{ 2/t,
	t/2, t, \infty \}$, where $t>2$.
\end{enumerate}

In the following subsection an example of how dilogarithm
 identities can be derived using the first feasible pair is given. Modifications and elaborations of this in subsequent sections gives the identities stated above.

\subsection{The identities}\label{ss:identities}
Let $(T, P)$ be a feasible pair given as in (i) of subsection \ref{ss:feasible}, $S$ the hyperbolic cylinder obtained from the identification of the two sides of $P$ by $T$.
Let $E$ be the collection of (unordered) pairs of nonadjacent sides of
the universal cover $\tilde
S$ of $S$. Then
\begin{equation}
E= \{ ( u, v)\,:\,  u, v \in \{ T^n(e_1), \, T^m (e_{\infty})\,
:\, n, m \in \Bbb Z\},\, d(u, v)>0\}, \end{equation}

\noindent
where $d(u, v)$ is the
hyperbolic distance
between $u$ and $v$, and $e_1= [1, t-1]$, $e_{\infty} = [t, \infty]$, where $[x_1,x_2]$ denotes the geodesic with endpoints $x_1,x_2$.
Note that $(u, v)$ and $(v, u)$ are considered as the same pair of geodesics.
The action of $T$ on $E$ splits $E$ into orbits. The following
is a set of representatives.
\begin{equation}
\mathcal E= \{( T^n(e_{\infty} ), e_{\infty})\,:\, n\ge 2\} \cup\{ (e_1, T^n (e_1))\,:\, n \ge 2\}
\cup \{(T^n(e_1), e_{\infty})\,:\, n\in \Bbb Z\}.\end{equation}
\noindent
Note that  the third set can be decomposed into
$ \{(T^n(e_1), e_{\infty})\,:\, n\ge 1\}
\cup  \{(e_1, T^n(e_{\infty}))\,:\, n\ge 1\}
\cup\{ (e_1, e_{\infty})\}$. This decomposition may
make the calculation easier on some occasions.
Applying Bridgeman's orthospectrum identity (see (3.2)), one has
\begin{equation}\sum_{(u,v)\in \mathcal E} \LL  ([ u, v]) = \pi^2/3,
\end{equation}

\noindent where $[u, v]$ is the cross ratio of $u$ and $v$ (note that here, $u$, $v$ are infinite geodesics and not points), defined from the cross ratio of their endpoints
(see Section \ref{s:crossratio}). The set
$\{ [u, v]\,:\, (u, v)\in \mathcal E\}$ is called {\em the set
	of cross ratios} of $(T, P)$, which can then be expressed in terms of known recurrence sequences.

\subsection{Outline of the paper}

Section 2 defines the cross ratio which will be used later and relates it to the distance between two non-intersecting complete geodesics. Section 3 gives Bridgeman's remarkable orthospectrum identity.
Section 4 onwards describes and proves the various infinite dilogarithm identities.

\medskip

\noindent {\it Acknowledgements}. We are grateful to Martin Bridgeman, Tengren Zhang and Sam Kim for their interest in this work and for helpful conversations and comments.

\section{Cross ratios and distances between complete geodesics }\label{s:crossratio}

Following the convention used in \cite{B2}, we define the cross ratio of $4$  points in $\hat{\mathbb{C}}$ (with at least three of them distinct) by
\begin{equation} [z_1,z_2,z_3,z_4]=\frac{(z_1-z_2)(z_4-z_3)}{(z_1-z_3)(z_4-z_2)}.\end{equation}
As is well known, the cross ratio is invariant under the action of elements of $PSL(2, \mathbb C)$.
If $x,y \in {\mathbb R}\cup \{\infty\}=\partial \mathbb H^2$, we will denote the geodesic from $x$ to $y$ by $[x,y]$. With this convention, if $x_1,x_2,x_3,x_4$ are four points in cyclic order in $\partial \mathbb H^2$ and $l$ is the perpendicular distance between the geodesics $[x_1,x_2]$ and $[x_3,x_4]$, then
\begin{equation}[x_1,x_2,x_3,x_4]=\frac{1}{\cosh ^2(l/2)}.
\end{equation}
If $u,v$ are  infinite geodesics with distance $l>0$ and end points $\{x_1,x_2\}$ and $\{x_3,x_4\}$  respectively such that $x_1,x_2,x_3,x_4$ are in cyclic order, then we define the cross ratio $[u,v]=[v,u]:=[x_1,x_2,x_3,x_4]$. If $u$ and $v$ have exactly one common endpoint (so distance  $l=0$), then $[u,v]=1$.

\section{Bridgeman's Identity }

\begin{figure}[h]
	\begin{center}
		\scalebox{0.3}[0.3]{\includegraphics{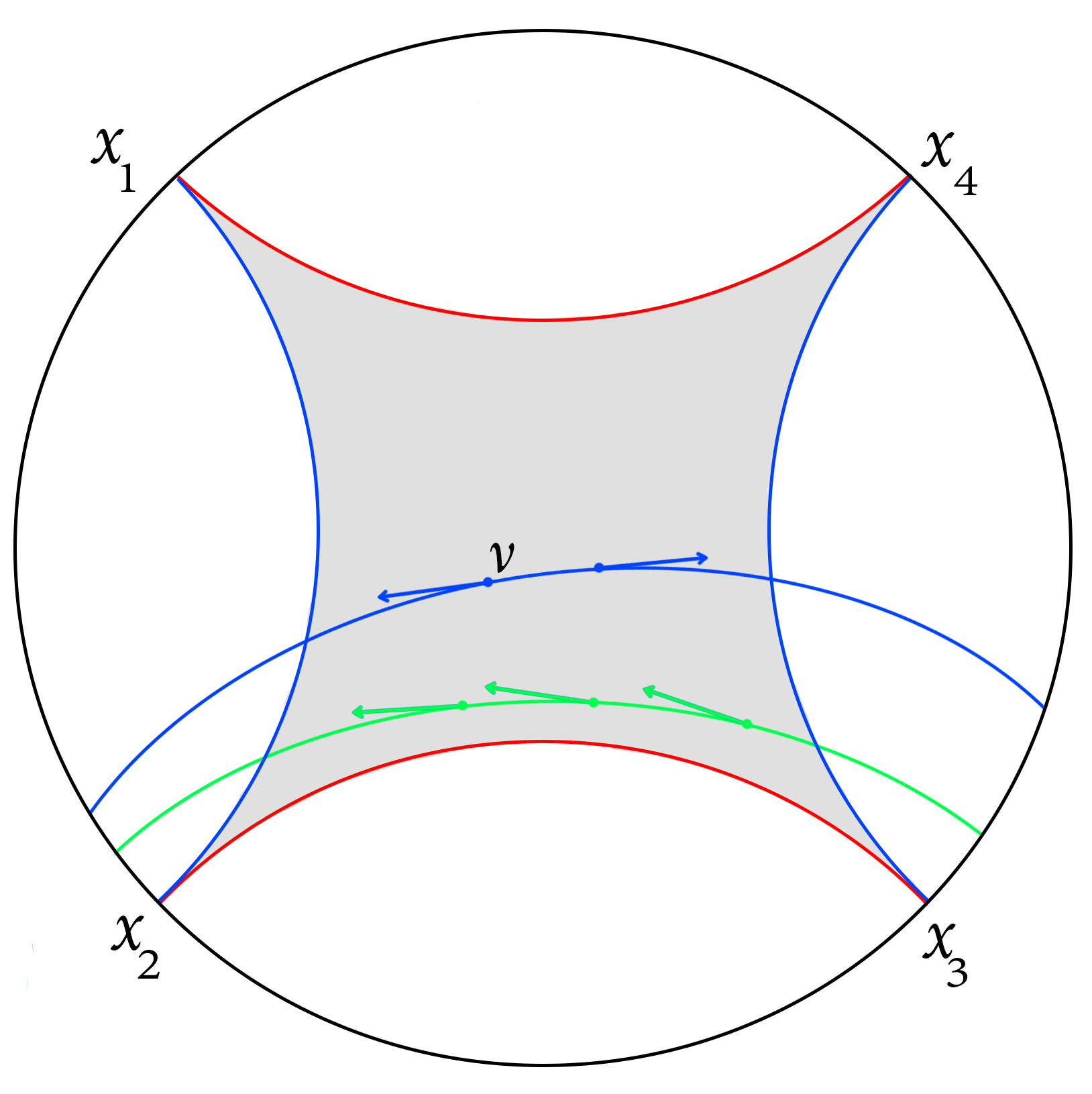}}
	\end{center}
	
	\caption{Unit vectors $v$ with base point in an ideal quadrilateral exponentiating to geodesics hitting opposite sides, using the disk model for $\mathbb H^2$.} \label{fig:Bridgeman}
\end{figure}

In \cite{B1}, Bridgeman showed that the measure of the  set $Q$ of unit tangent vectors whose base point lies in the ideal quadrilaterial with vertices $x_1,x_2,x_3,x_4 \in \partial \mathbb H^2$ (in cyclic order) and which exponentiate in both directions to a complete geodesic with one end point in the interval $(x_1,x_2) \subset \partial \mathbb H^2$ and the other endpoint in the interval $(x_3,x_4) \subset \partial \mathbb H^2$ is given by
\begin{equation}
\mu(Q)=8 \mathcal{L}([x_1,x_2,x_3,x_4]),\end{equation}
where $\mu(Q)$ is the measure of $Q$ and $\LL (x)$ is the Rogers dilogarithm defined at the beginning of the introduction, see figure \ref{fig:Bridgeman}.


More generally, he proved the following remarkable orthospectrum identity by decomposing the unit tangent bundle:

\begin{theorem} {\rm (Bridgeman \cite{B1})}
	For a finite area hyperbolic surface $S$ with totally geodesic boundary $\partial S \neq \emptyset$  and $N(S)$  boundary cusps, let $O(S)$ be the set of orthogeodesics in $S$, that is, the set of geodesic arcs with end points on $\partial S$, and perpendicular to $\partial S$ at both ends. Then
	\begin{equation}\sum_{\alpha \in O(S)} \LL  \left (\frac{1}{\cosh^2(l(\alpha)/2)}\right )=-\frac{\pi^2}{12}(6 \chi(S)+N(S))\end{equation}
	
	where  $\chi(S)$, the Euler characteristic of $S$,  satisfies $\chi(S)=-{\mathrm{Area}}(S)/2\pi$.
\end{theorem}

Generalizations of this in various settings  be found in \cite{BK}, \cite{BT} and \cite{LT}.

\section{The first feasible pair : (i) of subsection \ref{ss:feasible}}\label{s:firstpair}

\subsection{The first feasible pair}
We start with a basic case.
Let   $ \,{\tiny T =
\left ( \begin{array}{cc}
t&  -1\\
1 & 0\\
\end{array}
\right )}\in SL(2, \Bbb R)$, $(t>2$) be an element
of infinite order.
The two sides $[\infty,1]$ and $[t-1,t]$
of
the hyperbolic convex hull $P$ of $\{\infty, 1, t-1, t\}$ are
identified by $T$. Note that $T$ sends the interior of $P$ to the exterior of $P$. The identification of the two sides of $P$ by $T$ gives a hyperbolic cylinder $S$, a double crown with one boundary cusp on each boundary component.
The universal cover $\tilde S$ of $S$ is the hyperbolic
convex hull of $\{ T^n(1), T^n(\infty)\,:\, n \in \Bbb Z\}$.
The action of $T$ on the sides of $\tilde S$ gives two orbits
$\{T^n([t, \infty])\,:\,
n\in \Bbb Z\}$ and
$\{T^n([1, t-1])\,:\,
n\in \Bbb Z\}$.
Applying our discussion in subsection \ref{ss:identities},
the  set of cross ratios of pairs of non-adjacent sides of $\tilde{S}$ (see (1.2) and (1.3)) consists of the
following:

\begin{enumerate}
\item[(i)] $[(T^n(t), T^n(\infty), t, \infty]=[\infty, t, T^n(\infty), T^n(t)]=[\infty, T(\infty), T^n(\infty), T^{n+1}(\infty)]$,
where $n \ge 2$,
\item[(ii)] $[(1,t-1,  T^{n}(1), T^n( t-1)]=[1, T(1), T^n(1), T^{n+1}(1)]$, where $n \ge 2$,

\item[(iii)] $[(1, t-1 , T^n(t), T^n(\infty) ]=[(1,T(1) , T^{n+1}(\infty), T^n(\infty) ]$,
where $n \ge 1$,
\item[(iv)] $[(T^n(1), T^n(t-1) , t, \infty]=[T^n(1), T^{n+1}(1),T(\infty), \infty]$,
where $n \ge 1$,
\item[(v)] $[1, t-1, t, \infty ]=[1,T(1),T(\infty), \infty]$.

\end{enumerate}

\smallskip
\noindent
We shall now give a detailed study of $T^n$ so that one can
determine the cross ratios  above.
One can easily prove by induction that
\begin{equation}  T^n=
\left ( \begin{array}{cc}
q_n& -q_{n-1}\\
q_{n-1} & -q_{n-2}\\
\end{array}
\right ) , \quad  T^n(\infty) = {q_n\over q_{n-1}}, \quad T^n(1)= {(q_n-q_{n-1})\over (q_{n-1}-q_{n-2})}\end{equation}

\smallskip  \noindent
where
$  q_0=1$, $q_{1}=t$,
$q_n= tq_{n-1}-q_{n-2}, \quad n \in \mathbb Z.$

\smallskip

To simplify the calculations, we set
$p_n= q_n-q_{n-1}. $ It follows that $\{p_n\}$ is
defined by
$p_0=1, p_1=t-1, p_n=tp_{n-1}-p_{n-2}.$
It is then an easy matter to show that
\begin{equation}
p_kp_{k-2}=p_{k-1}^2+(t-2),\qquad \,p_n-p_{n-1}= (t-2)q_{n-1}.
\end{equation}

\noindent  We are now ready to calculate the cross ratios. The
cross ratios of (i)-(v) can be calculated easily
by the above mentioned recurrences. In particular, by direct calculation, we have
\begin{equation}
[1, t-1, t,\infty]=  [T(1), T(t-1), t , \infty].\end{equation}

\noindent Also, by (4.1) and (4.2), (i) and (ii) gives the same cross ratios. Indeed, we can always find a transformation $S\in PSL(2,\mathbb C)$ with the same fixed points as $T$ such that $S(\infty)=1$ and $TS=ST$  so the cross ratios in (i) and (ii) are the same by the invariance of cross ratios. Equivalently,
$\{ c \in \mbox{(i)}\} = \{ c\in
\mbox{(ii)}\}$. Furthermore,
\begin{equation} \sum _{c\in \mbox{(ii)}}\LL( c)
=\sum_{c\in \mbox{(i)}} \LL(c)
= \sum_{n=2}^{\infty }
\LL \left (
[ T^{n+1}(\infty), T^{n}(\infty), t, \infty]
\right) =
\sum_{n=2}^{\infty }
\LL \left (\left ( \frac {1}{q_{n-1}}\right ) ^2
\right)
,\end{equation}

\noindent where the first summation  represents the sum of all the
cross ratios coming from (ii).
Applying (4.1), (4.2) and (4.3) one has $\{ c \in \mbox{(iv)}\}
=\{ c \in \mbox{(iii)}  \cup \mbox{(v)  }\}$
and
\begin{equation}\sum_{c\in \mbox{(iv)}}\LL(c) =
\sum_{n=1}^{\infty }
\LL \left (
[ T^n(1), T^{n}(t-1), t, \infty]
\right)
=\sum_{n=0}^{\infty }
\LL \left ( \left (\frac{t-2}{(q_{n+1}-q_n)(q_{n-1}-q_{n-2})}
\right)\right )
.\end{equation}

\noindent Since all the set of cross ratios have been completely determined, we have:

\smallskip
\begin{proposition} Suppose that $t>2$. Then
\begin{equation}
\sum_{n=2}^{\infty }
\LL \left (\left ( \frac {1}{q_{n-1}}\right ) ^2
\right)  +
\sum_{n=0}^{\infty }
\LL \left ( \left (\frac{t-2}{(q_{n+1}-q_n)(q_{n-1}-q_{n-2})}
\right)\right )
= \pi^2/6,\end{equation}
\noindent where $\{q_n\}$ is the recurrence defined by
$  q_0=1 , q_1=t,
q_n= tq_{n-1}-q_{n-2}.$
\end{proposition}

\begin{proof} The pair $(T,P)$ has two boundary cusps and
area $2\pi$. By (4.4), (4.5) and  (1.3),
one has
$$
2 \sum_{n=2}^{\infty }
\LL \left (\left ( \frac {1}{q_{n-1}}\right ) ^2
\right)  +
2\sum_{n=0}^{\infty }
\LL \left ( \left (\frac{t-2}{(q_{n+1}-q_n)(q_{n-1}-q_{n-2})}
\right)\right )
= \pi^2/3.$$

 \end{proof}

\subsection{Fibonacci numbers}
Let  $ T= {\tiny
\left ( \begin{array}{cc}
3&  -1\\
1 & 0\\
\end{array}
\right )}$. Then
$
T^n =
{\tiny \left ( \begin{array}{cc}
f_{2n+2}& - f_{2n}\\
f_{2n} & -f_{2n-2}\\
\end{array}
\right )},$
where
$f_n$ is the $n$-th Fibonacci number.
Consequently, $T^n(\infty) = f_{2n+2}/f_{2n}$ and
$T^n(1) = f_{2n+1}/f_{2n-1}$. By (4.6), one has
\begin{equation}\sum_{k=1}^{\infty }
\left ( \LL \left (\frac{1}{f_{2k+2}^2}\right )
+ \LL \left (\frac{1}{f_{2k-3}f_{2k+1}}\right )\right ) = \frac{\pi^2}{6},\end{equation}

\subsection{Lucas numbers} The Lucas numbers are defined by
$l_0 = 2, l_1=1, l_n= l_{n-1}+l_{n-2}$. Let $ {\tiny T =
	\left ( \begin{array}{cc}
	7&  -1\\
	1 & 0\\
	\end{array}
	\right )}$
and $P$  the hyperbolic convex hull of
$\{\infty, 1, 6, 7\}$. Note that
$t=7$, $f_4=3$  and $f_4^2(t-2)=45$. By (4.6), and expressing the terms $q_n$ in terms of the Fibonacci and Lucas numbers, we have:

\begin{equation}\sum_{k=2}^{\infty }
\LL \left (\frac{3^2}{f_{4k}^2}\right)
+
\sum_{k=0}^{\infty }
\LL \left (\frac{45}{l_{4k-2}l_{4k+6}}\right)
= \frac{\pi^2}{6}.\end{equation}

\section {Decomposition of Proposition 4.1 }

The main purpose of this section is to
break (4.6)  of Proposition 4.1 into two dilogarithm identities, and to show that one of the resulting identities is equivalent to that obtained by Bridgeman in \cite{B2}. Geometrically what we do is to decompose the double crown into two crowns using the unique closed geodesic (``waist'') embedded in the double crown.

Let $(T, P)$ be the feasible pair  given as in subsection 4.1,
where
$T={\tiny
\left ( \begin{array}{cc}
t&  -1\\
1 & 0\\
\end{array}
\right )}$ and $P$ is the hyperbolic convex hull of
$\{\infty, 1, t-1, t\}$. We shall decompose $P$ as follows.
Let $ u  > 1/u $ be the fixed points of $T$.
Then $[1/u  , u]$ intersects
$[1, \infty] $ and $[t-1, t]$ at two
points $v_1$ and $v_2$, where $[1/u, u]\cap [1,\infty]=
\{v_1\}$. See figure \ref{fig:FP}.
Let $P_1$ be the hyperbolic convex hull of
$\{\infty, v_1, v_2, t\}$.
Then
$T$ pairs $[v_1,\infty]$ and $[v_2, t]$ and
$(T, P_1)$ is also a feasible pair. The side pairing of $P_1$ gives a crown $S_1$ with one tine. We may also consider the other polygon $P_2$ which is the convex hull of $\{1, t-1, v_2, v_1\}$ and the corresponding pairs $(T, P_2)$ and surface $S_2$. As the surfaces $S_1$ and $S_2$ are isometric, the identities arising from both are essentially the same.
\smallskip

Let $\tilde {S_1}$
be the universal cover of $S_1$. It follows that
the sides of $\tilde {S_1}$ has two orbits:
$\{[1/u , u]\}$ and
$\{ T^n( [t, \infty])\,:\, n\in \mathbb Z\}$.
Note that $T(\infty) = t=u+1/u>2$. By (3.2), one has
\begin{equation}
 \LL ([1/u , u, t, \infty]) +
\sum_{n=2}^{\infty }
\LL \left (
[ T^{n+1}(\infty), T^{n}(\infty), t, \infty]
\right)=\pi^2/6.\end{equation}

\noindent  Applying (4.4) and the fact that $\LL (x) +\LL(1-x) =\pi^2/6$, we have
\begin{equation}
\sum_{n=2}^{\infty }
\LL \left (\left ( \frac {1}{q_{n-1}}\right ) ^2
\right)  =\LL (1/u^2),\end{equation}
\noindent where $q_n$ is the  recurrence defined by
$q_0=1 , q_1=t, q_n= tq_{n-1}-q_{n-2}.$

\smallskip

\begin{theorem} Suppose that $t>2$ $($not necessarily
an integer$)$.
Let $u > 1/u$ be the roots
of $x^2-tx+1=0$. Then $t= u+ 1/u$ and
\begin{equation}
\sum_{n=1}^{\infty }
\LL \left ( \frac {1}{q_{n}\,^2}\right )
=\LL (1/u^2), \qquad
\sum_{n=1}^{\infty }
\LL  \left (\frac{t-2}{(q_{n}-q_{n-1})(q_{n-2}-q_{n-3})}
\right)=\LL ( 1-1/u^2),
\end{equation}

 \noindent where $q_n$ is the  recurrence defined by
   $q_0=1 , q_1=t, q_n= tq_{n-1}-q_{n-2}.$
\end{theorem}

\begin{proof} Apply  (5.2) and (4.6).
\end{proof}

\subsection{An equivalent form of Bridgeman's identity}

Bridgeman studied the universal cover of a hyperbolic surface $S$
which is topologically an annulus with one boundary component being
a closed geodesic of length $L >0$ and the other an infinite geodesic
with a single boundary cusp (a crown with one tine) and proved that
\begin{equation}\LL ( e^{-L})
= \sum_{k=2}^{\infty}
\LL  \left (\frac{\sinh ^2(L/2)}
{\sinh ^2 (kL/2)}\right ). \end{equation}

\noindent   In our notation, he used $ T= {\tiny
	\left ( \begin{array}{cc}
	\sqrt \lambda&  0\\
	0&1/\sqrt \lambda \\
	\end{array}
	\right )}$, $P$ the hyperbolic convex hull of $\{i,1,\lambda, \lambda i\}$ where $\lambda=e^L$, and $T$ identifies $[i,1]$ with $[\lambda i, \lambda]$ to get $S$. The universal cover $\tilde S$ under his study is the hyperbolic
convex hull of
$\{0, \infty\} \cup \{  T^{k-1}(\lambda)=\lambda^k\,:\, k \in \ZZ \}$,
since $T (x) = \lambda x.$
We show in this section
that (5.4) and (5.2) are equivalent, which is not surprising as the crowns obtained from the two constructions are isomorphic for appropriate choices of $L$ and $t$.

Adapting to Bridgeman's pair $(T,P)$ (his $\sqrt \lambda$
is our $u$), we set

\smallskip
$ t := \sqrt \lambda + 1/\sqrt \lambda >2$,
$ A:= {\tiny
\left ( \begin{array}{cc}
t&  -1\\
1 & 0\\
\end{array}
\right )}.$
By (4.1), one has
\begin{equation} A^{n}= {
\left ( \begin{array}{cc}
q_n& -q_{n-1} \\
q_{n-1} & -q_{n-2}\\
\end{array}
\right )},\end{equation}

\noindent where
$q_{-1}=0, q_{0}=1,  q_n= tq_{n-1}- q_{n-2}, \, n \in \ZZ.$
Note that $A$ is similar to $T$. We have:

\begin{proposition}
	Identities 
	$(5.4)$ and $(5.2)$ are equivalent.
	
\end{proposition} 

\begin{proof}
To simplify the calculations, we define a recurrence $v_n$ by
$$v_{-1} = 2, v_0=t,  v_n= tv_{n-1}-v_{n-2}, \quad n>1.$$
Let $T$, $P$  and $\tilde S$ be given as above and 
let $ V= {\tiny
\left ( \begin{array}{cc}
\sqrt \lambda &  1/\sqrt \lambda\\
1 & 1\\
\end{array}
\right )}.$  Then 
$$V^{-1} AV={\tiny
\left ( \begin{array}{cc}
\sqrt \lambda &  0\\
0 & 1/\sqrt \lambda\\
\end{array}
\right )}=T$$  and $V\tilde  S$ is  the  universal cover  of $VS$
that is invariant under the action of $A$.
The sides of $V\tilde S$ under the action of $A$ has two orbits
$\{ [1/\sqrt \lambda, \sqrt \lambda]\}$ and
$\{A^k V[1, \lambda]\,:\, k\in \ZZ \}$.
Applying
(1.3) and (3.4) to $V\tilde S$, (5.4)
is equivalent to
\begin{equation}
\LL (1/\lambda )
= \sum_{k=2}^{\infty}
\LL  \left (
[A^0V(1), AV(1), A^kV(1), A^{k+1}V(1)]
\right ),\,\, {where}\,\, V(1) = t/2.   \end{equation}

\noindent It is an easy matter to show that  $A^kV(1)= v_{k}/v_{k-1}$.
As a consequence, the above identity can be simplified as follows.
\begin{equation} \LL (1/\lambda )
= \sum_{n=2}^{\infty}
\LL  \left (\left ( \frac{ t^2-4}
{v_n-v_{n-2}}\right )^2\right ).\end{equation}

\noindent One can show by induction that $1/q_{n-1} = (t^2-4)/(v_n-v_{n-1})$.
This completes the proof of the proposition.
\end{proof}

\noindent {\it Remark:}  The difference between (5.4) and (5.2) is that $(5.4)$ arises from the
study of the Jordan canonical form
which express the functions in terms of the eigenvalues of
$T$ (see the proof of Theorem 2.1 of
\cite{B2}) while (5.2)
arises from the study of $T$ in cyclic basis. The surfaces are of course isometric so that the identities have to be identical, the difference is that one is expressed in terms of the hyperbolic sine of the length $L$ of the waist, the other is expressed in terms of recurrences defined from the trace $t$, where $t=2\cosh( L/2)$.

\subsection{Examples} We give some examples of how Theorem 5.1
works.

\begin{example} Let  $\phi = (1+\sqrt 5)/2 $ and let
$f_1=1,f_2=1, f_3=2, f_4 = 3
$ be the Fibonacci numbers. Then
\begin{equation}
\sum_{k=2}^{\infty}
\LL  \left (\left ( \frac{ 1}
{f_{2k}}\right )^2\right )=\LL (1/\phi^4 ), \qquad
\sum_{k=1}^{\infty }
\LL \left (\frac{1}{f_{2k-3}f_{2k+1}}\right ) =
\LL (1-1/\phi^4).\end{equation}

\end{example}

\noindent The above  is a special case
of our results in Section 8 (see (8.7)) and comes from our study of $u = \phi^2,$
$t=3$ in Theorem  5.1. In the  case $u = \phi^4$, one has $t=7$.
Theorem 5.1 gives
\begin{equation}
\sum_{k=2}^{\infty }
\LL \left (\frac{3^2}{f_{4k}^2}\right )=\LL (1/\phi^8)
 ,\qquad
\sum_{k=0}^{\infty }
\LL \left (\frac{45}{l_{4k-2} l_{4k+6}}\right)=\LL (1-1/\phi^8),
\end{equation}

\noindent where $l_n$ is the $n$-th Lucas number
($l_0=2, l_1=1, l_n= l_{n-1}+l_{n-2}$).

Formulas for $\LL (1/\phi^k)$ for general $k$ can be found in (8.7) and Sections 10-11.
The following gives $\LL (1/\phi^2)$.

\begin{example}  Let
$t =\phi
+1/\phi  = \sqrt 5$. By (5.2), one has

\begin{equation}\sum_{n=1}^{\infty}
\LL  \left (\left ( \frac{ 1}
{5f_{2n}}\right )^2\right )+
\sum_{n=2}^{\infty}
\LL  \left (\left ( \frac{ 1}
{l_{2n-1}}\right )^2\right )=\LL (1/\phi^2 )
.\end{equation}
\end{example}

\noindent
Identity (5.10) was first proved by Bridgeman [B2].
However, the current feasible pair  does not work well for
$\LL  (1-1/\phi^2)$ as $(t-2)/(q_n-q_{n-1})(q_{n-2}-q_{n-3})$ is not rational. See Section 12 for an identity for $\mathcal
L(1-1/\phi^2)$ expressed in terms of dilogarithms of rationals.

\section{Identities for two term recurrences
}

Let   $p_n$ be  the recurrence defined by
$p_{-2}=0$, $p_{-1}=1$, $p_n= ap_{n-1}+bp_{n-2}$ where $a$ and $b$ are
positive integers.
It follows that
\begin{equation}
\left ( \begin{array}{cc}
p_{n+2} \\
p_{n+1} \\
\end{array}
\right )
=
\left ( \begin{array}{cc}
a^2+b & ab\\
a & b\\
\end{array}
\right )
\left ( \begin{array}{cc}
p_n \\
p_{n-1} \\
\end{array}
\right ).
\end{equation}

\noindent  Let $T={\tiny  \left ( \begin{array}{cc}
(a^2+b )/b& a\\
a/b & 1\\
\end{array}
\right ) \in PSL(2,\mathbb R) }$. Then  $T$ is similar
to
${\tiny R=: \left ( \begin{array}{cc}
a^2/b  +2& -1\\
1 & 0\\
\end{array}
\right ) }$.
Let $P$ be the hyperbolic convex hull
of $\{\infty, 1, {a^2\over b}+1, {a^2 \over b}+2\}$. Then $(R.P)$ is a feasible pair and we call the dilogarithm identity
associated with
$(R, P)$ the dilogarithm identity of the recurrence
$p_n$.

\begin{proposition}
Let $a$ and $b$ be positive
integers and let $\{p_n\}$
be given as in $(6.1)$.
Then the dilogarithm identity associated with $\{p_n\}$ is
\begin{equation}
\LL (1/u^2 )
= \sum_{n=1}^{\infty}
\LL  \left ( \frac{1}
{q_n^2}\right ),
\end{equation}

\medskip

\noindent where $ u +1/u=t = a^2/b +2 $ is the trace of $T$, $u >1$,  and
$q_n$ is the recurrence defined by $q_0=1, q_1=t, q_n=tq_{n-1}-
q_{n-1}$. Furthermore, 
\begin{equation}
\frac{1}{q_n}=\frac{b^n}{p_{2n-1} +b p_{2n-3} + \cdots + b^k p_{2n-2k-1} + \cdots+ b^n p_{-1}}. \end{equation}
\end{proposition}

\begin{proof}
$T$ is similar to
${\tiny  \left ( \begin{array}{cc}
a^2/b +2& -1\\
1 & 0\\
\end{array}
\right ) }$.
The proposition can be proved by applying (5.2). The expression in terms of $p_n$ is a straightforward computation.
\end{proof}

\begin{example} Let $a=2$ and $b=3$. Then
the $p_n$'s (starting at $p_{-2}$) are given as follows.
$$
\begin{array}{ccccccccc}
p_n:& 0 &1&2&7&20&61&182 &\cdots\\
\end{array}
$$

\noindent It follows that
$T={\tiny  \left ( \begin{array}{cc}
7/3& 2\\
2/3 & 1\\
\end{array}
\right ) }$ and  $T$ is similar to ${\tiny  \left ( \begin{array}{cc}
10/3& -1\\
1 & 0\\
\end{array}
\right ) }$.
By Proposition 6.1, one has
\begin{equation}
\mathcal
L\left (\frac{1}{9}\right )
= \sum_{k=1}^{\infty} \LL \left (\left (
\frac {3^{k}}{ 9^k +9^{k-1}+\cdots +9^2+9+1}\right ) ^2\right ).\end{equation}
\end{example}

\noindent {\em Remark} : This identity can be interpreted as arising from the convergents of a non-standard continued fraction expansion arising from a non-arithmetic lattice $\Gamma(1,1/3,3)$, see \cite{LTV}, section 7 for the definition of $\Gamma(1,1/3,3)$ and more details. A more general version is given below.

\section {Identities for $\LL (1/n)$
and Chebyshev polynomials }
\subsection{Identities for $\LL (1/n)$}
Let $u =\sqrt n$, $t = \sqrt n +1/\sqrt n$, where $n >2$. (Of particular interest is when $n$ is an integer).
Identity  (5.2) gives
\begin{equation}
\mathcal
L\left (\frac{1}{n}\right )
= \sum_{k=1}^{\infty} \LL \left (\left (
\frac {n^{k/2}}{ n^k +n^{k-1}+\cdots +n^2+n+1}\right ) ^2\right ).\end{equation}

\subsection{Chebyshev polynomials }

Let $q_n$ be the recurrence defined in Theorem 5.1 and let $t=2x$.
It follows that
\begin{equation}
q_n = U_n(x),\end{equation}

\noindent
where $U_n(x)$  is the $n$-th Chebyshev polynomial
of the second kind.
Applying (5.3) gives, for $x>1$:
\begin{equation}
\begin{aligned}
\sum_{n=1}^{\infty }
\LL \left ( \frac {1}{U_n(x)\,^2}\right )
&=\LL \left(\frac{1}{(x+\sqrt{x^2+1})^2}\right), \qquad \quad\\
\sum_{n=1}^{\infty }
\LL  \left (\frac{2x-2}{(U_{n}(x)-U_{n-1}(x))(U_{n-2}(x)-U_{n-3}(x))}
\right)&=\LL \left( 1-\frac{1}{(x+\sqrt{x^2+1})^2}\right). \qquad 
\end{aligned}
\end{equation}
 The first identity of (7.3) was derived by Bridgeman in \cite{B2}.

\section {The second feasible pair : (ii) of subsection 1.3 }

The feasible pair  $(T, P)$ given as in Proposition 4.1 does not work well
for continued fractions
as one cannot tell if the $q_n$'s in Theorem 5.1 are related to the
$n$-th convergents of $u$. We will develop another feasible pair
that works well for Fibonacci numbers (subsection 8.2) as well as
continued fractions (Section 9).

\subsection{The second feasible pair}

Let
$A={\tiny
\left ( \begin{array}{cc}
a& c\\
b & d\\
\end{array}
\right ) }\in SL(2, \Bbb Z)$ be a matrix with positive entries and trace $a+d >2$.
$A$ identifies the sides $[0, \infty]$ and
$[c/d, a/b]$
of the hyperbolic convex hull $P$ of $\{0, c/d, a/b, \infty\}$
and  sends the interior of $P$ to the exterior of $P$ so $(T,P)$ is a feasible pair.
Let $t= a+d$ be the trace of $A$. Then $A^n= tA^{n-1}-A^{n-2}$
for all $n.$
Hence, if we set
$A^{n}=
{\tiny \left ( \begin{array}{cc}
p_{2n-1}& p_{2n-2}\\
q_{2n-1} & q_{2n-2}\\
\end{array}
\right )} $, where
$A^{0}=
{\tiny \left ( \begin{array}{cc}
p_{-1}& p_{-2}\\
q_{-1} & q_{-2}\\
\end{array}
\right ) =
\left ( \begin{array}{cc}
1& 0\\
0 & 1\\
\end{array}
\right )} ,$ then
\begin{equation}
A^n =  \left ( \begin{array}{cc}
p_{2n-1}& p_{2n-2}\\
q_{2n-1} & q_{2n-2}\\
\end{array}
\right )
= t
\left ( \begin{array}{cc}
p_{2n-3}& p_{2n-4}\\
q_{2n-3} & q_{2n-4}\\
\end{array}
\right )
-
\left ( \begin{array}{cc}
p_{2n-5}& p_{2n-6}\\
q_{2n-5} & q_{2n-6}\\
\end{array}
\right ).\end{equation}

\noindent This gives recurrence formulas for the $p_n$'s and $q_n$'s.
One can show by induction that
the $q_k$'s and $p_k$'s satisfy the following
identities for all $n$.
\begin{equation}
\det\, {
\left ( \begin{array}{cc}
p_{2n+1}& p_{2n+1+k}\\
q_{2n+1} & q_{2n+1+k}\\
\end{array}
\right ) } = q_{k-1},\,\,\det\, {
\left ( \begin{array}{cc}
p_{2n}& p_{2n+k}\\
q_{2n} & q_{2n+k}\\
\end{array}
\right ) } = -p_{k-2}.\end{equation}

\noindent
Note that  $p_0/q_0=c/d$, $p_1/q_1= a/b$.  Note also that
\begin{equation}
A^{-n} =
{\left ( \begin{array}{cc}
p_{-2n-1} & p_{-2n-2}\\
q_{-2n-1} & q_{-2n-2}\\
\end{array}
\right ) } = (A^n)^{-1} =
{
\left ( \begin{array}{cc}
q_{2n-2} &- p_{2n-2}\\
-q_{2n-1} & p_{2n-1}\\
\end{array}
\right ) }.\end{equation}

\noindent
The sides of the universal cover $\tilde S$ of the cylinder $S$
form two orbits $\{ A^n([a/b, \infty])\,:\, n \in \Bbb Z\}$
and $\{ A^n([0, c/d])\,:\, n \in \Bbb Z\}$. The following identity holds:

\begin{multline}
\sum _{n=2}^{\infty}
\LL \left (\left (\frac {b}{q_{2n-1}}\right)^2 \right )
+
\sum _{n=2}^{\infty}
\LL \left (\left (\frac {c}{p_{2n-2}}\right)^2 \right )
+\sum_{n=1}^{\infty }\LL \left (
\frac{bc}{q_{2n}q_{2n-4}}\right ) +
\sum_{n=1}^{\infty }\LL \left (
\frac{bc}{p_{2n+1}p_{2n-3}}\right )\\
+ \LL ([0, c/d, a/b ,\infty]) = \pi^2/3.
\end{multline}

\begin{proof} Similar to Proposition 4.1,
the  cross ratios associated with
$(T, P)$ consists of the following sets of cross ratios (see (1.3)).

\begin{enumerate}
	\item[(i)] the cross ratios $[(A^n(a/b), A^n(\infty), a/b, \infty]$
	and  $[(0,c/d,  A^{n}(0), A^n( c/d)]$,
	where $n \ge 2$,
	\item[(ii)]
	the cross ratios $[(0, c/d , A^n(a/b), A^n(\infty) ]$
	and $[(A^n(0), A^n(c/d) , a/b, \infty]$,
	where $n \ge 1$,
	\item[(iii)]
	the cross ratio $[0, c/d, a/b, \infty ]= bc/ad$.
	
\end{enumerate}

\noindent    Apply the above  and identity
(1.3), where the terms in the infinite sums correspond to (i) and (ii) above and are expressed in terms of the $p_n$'s and $q_n$'s. \end{proof}


Recall the way we decompose identity (4.6) into two identities
given as in Theorem 5.1. A similar decomposition can be done to (8.4).
Let $A$ be given as above and let $w >\overline w$ be
the fixed points of $A$. Let $v_1$ and $v_2$
be the intersection
$[0, \infty]\cap [\overline w, w]$ and
$ [c/d, a/b]\cap[\overline w, w]$ respectively.
$A$ identifies $[v_1, \infty]$ and $[v_2, a/b]$ and
$(A, P_1)$ is a feasible pair where
$P_1$ is the hyperbolic convex hull of
$\{\infty, v_1, v_2, a/b\}$. The surface $S_1$ obtained from $P_1$ is a crown with one tine.
Lift to the universal cover of $S_1$.
The sides of $\tilde S_1$ split into
two orbits $\{(\overline w, w)\}$
and $\{ A^n((a/b ,\infty))\,:\, n\in \Bbb Z\}$.
Similar to the way we get (5.2),
we apply (1.3) and (3.2)
to
get the following identity.
\smallskip

\begin{proposition} Let $b>0$ be an integer
and let $A$ and $q_n$ be given as in $(8.1)$.  Then

\begin{equation}
\LL (1/u^2) =          \sum _{n=2}^{\infty}
\LL \left (\left (\frac {b}{q_{2n-1}}\right)^2 \right )
\end{equation}

\smallskip
\noindent where $u> 1/u$ are the eigenvalues of $A$.
\end{proposition}

\smallskip
We find $(8.5)$  works better for $\LL (1/\phi^{4k})$
(see Remark $8.3)$  and continued fractions $($see Section 9$)$.

\begin{remark}
Note   that the  cross
ratios below satisfy the following interesting properties.
\begin{enumerate}
\item[(i)]
$\mathcal [A(0), A^{2}(0), a/b, \infty]
=[0, c/d, a/b, \infty] = bc/ad$ if $d=1$,
\item[(ii)]
$[0, c/d, a/b, \infty]
=[0, c/d, A^2(\infty) , A(\infty)] =bc/ad$,
if $a=1$.

\end{enumerate}
Also, since we assume $a+d>2$, one cannot have $a=d=1$. Note also that
  $ A^n =  {\tiny
\left ( \begin{array}{cc}
 ar_n-r_{n-1} &  cr_n\\
 br_n & dr_n-r_{n-1}\\
\end{array}
\right )}$, where $r_0 =1, r_1=t = a+d, r_n =tr_{n-1}-r_{n-2}$.
Hence $bp_{2n-2} = cq_{2n-1}$.
 It follows that (8.4) can be expressed in the following form.
\begin{equation}
2\sum _{n=2}^{\infty}
\LL \left (\left (\frac {b}{q_{2n-1}}\right)^2 \right )
+\sum_{n=1}^{\infty }\LL \left (
\frac{bc}{q_{2n}q_{2n-4}}\right ) +
\sum_{n=1}^{\infty }\LL \left (
\frac{bc}{p_{2n+1}p_{2n-3}}\right )
+\LL \left (\frac{bc}{ad}\right ) = \pi^2/3.\end{equation}

\end{remark}

\subsection{Dilogarithm identity for $\LL (1/\phi^{4k})$}

Identity (5.8) studies the matrix
$  {\tiny
\left ( \begin{array}{cc}
3&  -1\\
1 & 0\\
\end{array}
\right )}$ which is a conjugate of
${\tiny \left ( \begin{array}{cc}
f_3&  f_2\\
f_2 & f_1\\
\end{array}
\right )
}$ whose trace is 3 and gives
an identity for $\mathcal (1/\phi^4)$.
It is therefore  natural  to extend our
study to
$  {\tiny \left ( \begin{array}{cc}
a&  c\\
b & d\\
\end{array}
\right )=
\left ( \begin{array}{cc}
f_{2k+1}&  f_{2k}\\
f_{2k} & f_{2k-1}\\
\end{array}
\right )
}$.  By Proposition 8.1, one has,
\begin{equation}
 \LL \left( \left (\frac{1}{\phi}\right )^{4k}\right )
=\sum_{n=2}^{\infty }
\LL \left (\left (
\frac{f_{2k}}{f_{2nk}}\right )^2\right ).\end{equation}

\noindent Note that $b = f_{2k} $ and that the arguments of $\LL $ in (8.7) are all rational numbers with reduced forms $1/x$ for some $x \in \Bbb N$. However, Identity  (8.5)  does not work for $\LL (1/\phi^ {2(2k+1)})$
and $\LL (1/\phi^ {2k+1})$.
See Sections 10 and 11 for identities for $\LL (1/\phi^ {2(2k+1)})$
and $\LL (1/\phi^ {2k+1})$.

\begin{remark}  One can  also obtain Identity (8.7)
by applying Theorem 5.1. However, one  only gets
$\{q_i\} = \{ 1, t, t^2-1 ,\cdots \}$ (see Theorem 5.1 for the recurrence for $q_i$) rather than
$ \{f_{2k}/f_{2nk}\}=\{ f_{2k}/f_{4k}, f_{2k}/f_{6k}, f_{2k}/f_{8k},\cdots\}$
which we feel gives more insight and a direct connection with the Fibonacci numbers.
\end{remark}

\section{Dilogarithm identities for Continued fractions of period two }

\subsection{Periodic continued fractions} Let
$a_i$ $(i\ge 0$) be
 positive integers and let
\begin{equation}
\alpha =
    [{ a_0, a_1, a_2, \cdots, a_{l-1}}, \cdots ]
   = a_0+  \frac{1}
 {a_1 +\frac{1}{a_2 + \cdots}}\end{equation}

\smallskip
\noindent be a continued fraction.
 Let $r_n= p_n/ q_n
 =[
{ a_0, a_1, a_2, \cdots, a_{r}} ]
 $ ($n\ge 0$) be the $n$-th convergent of $\alpha$.
 Set $(p_{-2},  p_{-1})=(0,1) $
 and $(q_{-2}, q_{-1})= (1,0)$. Then
    $p_n$ and $q_n$ can be calculated recursively by
  $p_n= a_n p_{n-1}+ p_{n-2}$ and $q_n= a_n q_{n-1}+ q_{n-2}$.
If the $a_i$'s  satisfies $a_{n+l} = a_n$ for all $n \ge0$, we say that
 $\alpha$ is periodic of period $l$ and write
$\alpha =
    [  \overline
{ a_0, a_1, a_2, \cdots, a_{l-1}} ]$.

\subsection{Continued fractions of period 2}
Let
 $\alpha=[\overline {a,b}]$ be a continued fraction of
period 2 and let $r_n = p_n/q_n$ be the $n$-th convergent of $\alpha$.
Set
$A=  {\tiny
\left ( \begin{array}{cc}
ab+1&  a\\
b &1\\
\end{array}
\right )}
=
{\tiny
\left ( \begin{array}{cc}
p_1&  p_0\\
q_1&q_0\\
\end{array}
\right )
}.$
One can prove by induction that
$A^{n}= {\tiny  {
\left ( \begin{array}{cc}
p_{2n-1}&  p_{2n-2}\\
q_{2n-1} &q_{2n-2}\\
\end{array}
\right )}}.$
Let   $u >1/u$ be  the  eigenvalues of $A$. Similar to Proposition 8.1, we have:
\begin{equation}
 \LL (1/(b\alpha+1)^2) =
\LL (1/u^2) =          \sum _{n=2}^{\infty}
\LL \left (\left (\frac {b}{q_{2n-1}}\right)^2 \right ).
\end{equation}

\noindent Note  that the arguments of $\mathcal
L$ in the infinite sum involve the $(2n-1)$-th convergents of the continued fraction of $\alpha$ which gives
a close connection between the dilogarithm identity
associated with $u$ and the continued fraction of $\alpha$. When $b=1$, $u=\alpha+1$ is the solution of a positive Pell's equation and the identity was derived by Bridgeman in \cite{B2}.
Applying (8.6), gives a dilogarithm identity for $\LL  (
1-1/u^2)$, details are left to the reader.

%
%
%
%

\subsection{Generalised continued fractions}
Let $t >2$.
Define   
\begin{equation} \left <\overline  t\right >=
t+  \frac{-1}
{t +\frac{-1}{t + \cdots}}\end{equation}

\noindent $\left <\overline t\right > $ is called {\em the generalised continued fraction}
associated with $t$.
Define similarly the generalised $n$-th
convergent of $\left < \overline t\right >$. Let $u$ and $t$ be given as in
Theorem 5.1. Then $u =\left <\overline t\right > $  and
the generalised $n$-th convergent
of $\left <\overline t\right > $ is $q_n/q_{n-1}$.
As a consequence, Theorem 5.1
gives two   dilogarithm identities
(one for $\LL (1/u^2)$, one for
$\LL (1-1/u^2)$) where
the arguments of $\LL  $ of the infinite part are in terms of the
generalised $n$-th convergents of $u$.

\section{Identities for $\LL  (1/\phi^{4k+2})$
}
Let  $\phi= (1+\sqrt 5)/2.$
Since the
matrix $A= {\tiny\left ( \begin{array}{cc}
f_{2n}&  f_{2n-1}\\
f_{2n-1} &f_{2n-2}\\
\end{array}
\right )}$  has determinant $-1$ (the $A$ in  subsection 8.2
has determinant 1),
one cannot directly apply
(8.5)
to get  identities for $\LL (1/\phi^{4k+2})$ and $L(1/\phi^{2k+1})$.
However, identities for
$\LL  (1/\phi^{4k+2})$
can be obtained easily by
applying Theorem 5.1 as follows.
Let $u = \phi^{2k+1} $. Then $t =
u+1/u= f_{2k+1} \sqrt 5$. One can easily show by induction that the $q_n$'s in Theorem 5.1 with
$t = f_{2k+1} \sqrt 5$ satisfies
$q_{2(n-1)} = l_{n(4k+2)-(2k+1)}/l_{2k+1}$ and
$q_{2n-1} = \sqrt 5f_{n(4k+2)}/l_{2k+1}$. By (5.2), we get:

\begin{equation}
\LL  (1/\phi^{4k+2})
= \sum_{n=2}^{\infty} \LL  \left ( \frac{l_{2k+1}^2}{l_{n(4k+2)-(2k+1)}^2}\right )
+
\sum_{n=1}^{\infty} \LL  \left ( \frac{l_{2k+1}^2}{5f_{n(4k+2)}^2}\right ).
\end{equation}

\smallskip
\noindent
Note that the arguments of $\LL $ in the right hand side of (10.1)
are  rational numbers with reduced forms $1/x$ for some $x\in \Bbb N$.

\section{Identities for  $\LL  (1/\phi^{2k+1})$  }\label{s:twokplusone}

Let $u = \phi^{(2k+1)/2}$ and
$t = u+1/u$. Then
$t^2 -2= f_{2k+1}\sqrt 5$.
Similar to (10.1),  the $q_n$'s in Theorem 5.1 can be
split into two recurrences $H_n$ and $K_n$. To be more precise,
one has
\begin{equation}\LL  (1/\phi^{2k+1})
= \sum_{n=1}^{\infty} \LL  \left ( \frac{1}{
t^2H_n^2}\right )
+
\sum_{n=1}^{\infty} \LL  \left ( \frac{1
}{K_n^2}\right ),
\end{equation}

\noindent where
$t_0= t^2-2=f_{2k+1}\sqrt 5$, and
$H_n$ and $K_n$ are recurrences defined by
the following:
\begin{equation}
H_1=1,\,\, H_2=  t_0,  \,\,
H_n= t_0H_{n-1} - H_{n-2},\end{equation}
\begin{equation}
K_0=1,\,\, K_1= t_0 +1 ,\,\,
K_n= t_0K_{n-1}-K_{n-2}.\end{equation}

\smallskip
\noindent Similar to the way we split $\{q_n\}$ into
$\{H_n\}$ and $\{K_n\}$,   $\{H_n\}$ can be split
further. To be more precise,
the first summation
of (11.1) splits into
\begin{equation}
%
\sum_{n=1}^{\infty} \LL  \left ( \frac{1}{
t^2H_n^2}\right )
   =
\sum _{n=1}^{\infty} \LL  \left ( \left (
\frac{l_{k_0}}{ t\cdot l_{(2n-1)k_0} } \right ) ^2 \right )
+
\sum _{n=1}^{\infty} \LL  \left ( \left (
\frac{l_{k_0}}{\sqrt 5 t\cdot f_{2nk_0} } \right ) ^2 \right ),
\end{equation}

\noindent where $t^2 = f_{2k+1} \sqrt 5+2$, $k_0 = 2k+1$. This implies that
the arguments of $\LL $ are not rational. This is very different
from (10.1) and (8.7) where the arguments of $\mathcal  L$ in
those identities are rational numbers.
Returning to (11.1), similar to the first summation,  one can show that
the second summation of (11.1)  splits into

{\small \begin{equation}\sum_{n=1}^{\infty} \LL  \left ( 1/K_n^2\right )=
\sum _{n=1}^{\infty} \LL  \left ( \left (
\frac{l_{k_0}}{ l_{(2n-1)k_0} +\sqrt 5f_{2nk_0} } \right ) ^2 \right )+
\sum _{n=2}^{\infty} \LL  \left ( \left (
\frac{l_{k_0}}{l_{(2n-1)k_0} +\sqrt 5
	f_{2(n-1)k_0} } \right ) ^2 \right )
\end{equation}}

\smallskip
\noindent
where as before, $k_0 = 2k+1$.
This completes our study of the dilogarithm identities for
$\LL  (1/\phi^{2k+1})$.

\section{The third feasible pair : (iii) of subsection 1.3}

Consider the  feasible pair $(T, P)$, where $ {\tiny T =
\left ( \begin{array}{cc}
t&  -1\\
1 & 0\\
\end{array}
\right )}$ and $P$ is the hyperbolic convex hull of $\{\infty,
2/t, t/2, t\}$. Similar
to Proposition 4.1, one has
\begin{equation}
\sum_{n=1}^{\infty}
\left (
2\LL  \left( \left (\frac{1}{q_n}\right )^2 \right )+
\LL   \left (\frac{t^2-4}{p_np_{n-2}}\right)  +
\LL  \left (\frac{t^2-4}{p_{n+1}p_{n-1}}\right )
\right )
=\pi^2/3,\end{equation}

\noindent where
$p_0=2, p_1 = t,  p_n=  t p_{n-1}-p_{n-2}
\mbox{  and  }
q_0=1, q_1 =t,  q_n= t q_{n-1}-q_{n-2}.$ Note
that $q_{n}-q_{n-2}= p_n$ and that $ p_{n+2}-p_n= (t^2-4)q_n$.
Note also that $T^n(\infty) = q_n/q_{n-1},$  $ T^n(2/t) =p_n/p_{n-1}$.

\smallskip
\noindent {\bf Example 12.1.} In the case $t=\sqrt 5$, identity (12.1)
reduces to  the following identity.
\begin{equation}
\sum_{k=1}^{\infty } \left (
\LL \left (\frac{1}{5f_{2k}^2}\right)
+
\LL \left (\frac{1}{l_{2k+1}^2}\right)
+
\LL \left (\frac{1}{l_{2k-2} l_{2k}}\right)
+
\epsilon_k   \LL \left (\frac{1}{5f_{2k-3}f_{2k-1}}\right)\right )
= \pi^2/6
,\end{equation}

\noindent where $\epsilon _1= 1/2$ and $\epsilon _k= 1$
otherwise. Comparing to $(5.10)$ and  $(12.2)$ gives
us  an identity for $\LL(1-\/\phi^2)=\LL(1/\phi)$ in terms of dilogarithms of rational numbers which is not accessible by the previous methods. We have:
\begin{equation}
\sum_{k=1}^{\infty }\left (
\LL \left (\frac{1}{l_{2k-2} l_{2k}}\right)
+
\epsilon_k \LL \left (\frac{1}{5f_{2k-3}f_{2k-1}}\right)\right )
 =\LL (1-1/\phi^2)=\LL (1/\phi)=\frac{\pi^2}{10}.\end{equation}

\section{Identities associated  with
	$T \,:\, x\to x+1$  }

For this section it is convenient to introduce orthogeodesics of length $0$ between adjacent infinite geodesics meeting at a boundary cusp, that is, one orthogeodesic of length $0$ for each boundary cusp. This will allow for more compact expressions for the identities obtained. We define the enlarged set $\hat{O}(S)$ of orthogeodesics to be the set of orthogeodesics with these orthogeodesics of zero length included. Then, since $\LL(1)=\pi^2/6$, we can re-write Bridgeman's identity (3.2) as
\begin{equation}
\sum_{\alpha \in \hat{O}(S)} \mathcal L \left(\frac{1}{\cosh^2(l(\alpha)/2)}\right)=-\frac{\pi^2}{12}(6 \chi(S)-N(S)).
\end{equation}

\medskip

Consider the ideal hyperbolic polygon $P$ with  vertices $v_1,\ldots, v_{n}$, $n \ge 3$ with  $ v_1=0/1 < v_2 <\ldots< v_{n-1}=1/1$, $v_n=\infty$.
Identifying the two vertical sides of $P$ by the translation $T(z)=z+1$ gives a hyperbolic cylinder $S$ with fundamental
domain $P$ so $(T,P)$ is a feasible pair.
Note that $S$ is a degenerate crown, it has one regular cusp and $n-2$ boundary cusps.
Lifting $S$ to the universal cover, we see that
the enlarged set of orthogeodesics $\hat O(S)$ in (13.1)
corresponds to pairs of geodesics of the
following form:

\begin
{enumerate}
\item[(i)] $[v_i, v_{i+1} ]$ and $[ v_j , v_{j+1} ]$, $1\le j < i  \le n-2,$
\item[(ii)] $[v_i, v_{i+1}]$ and $[k+ v_j , k+v_{j+1}]
=T^k[v_i, v_{i+1}] $, $1\le i, j \le n-2,$
$k \ge 1$.
\end{enumerate}

%
%
%


Since $S$ has $n-2$ boundary cusps and $\chi(S) = (2-n)/2$,
 identity (13.1) gives

\begin{multline}
\mathop{\sum\sum}_{1 \le j<i\le n-2}
\LL \left(
[v_i,v_{i+1},v_{j},v_{j+1}]
\right)+ \sum_{i=1}^{n-2}\sum_{j=1}^{n-2} \sum_{k=1}^{\infty}
\LL \left(
[v_i,v_{i+1},k+v_{j},k+v_{j+1}]
\right)
=\frac{(n-2)\pi^2}{3}. 
\end{multline}

%
%
%
%
%
%

\medskip
\noindent {\bf Example 13.1.}
Let $n=3$. Then $v_1=0,$  $v_2=1$, $v_3=\infty$, there are no terms in the double sum and the triple sum reduces to a single sum, so (13.2) gives the  Richmond Szekeres identity (see also \cite{McS} where the identity was derived in the same way as here).
\begin{equation}
\sum_{k=2}^{\infty} \LL  \left (\frac{1}{k^2}\right ) = \pi^2/6.\end{equation}

Let $n=4$. The more interesting cases are
$v_1=0,$  $v_2 = f_{k-1}/f_k$, $v_3=1$, $v_4=\infty$ and
$v_1=0,$  $v_2 = 1/\phi$, $v_3=1$, $v_4=\infty$,
where $f_n$ is the $n$-th Fibonacci number and
$\phi = (1+\sqrt5)/2$ is the golden ratio. Letting $k \rightarrow \infty$ in the first case gives the second case, which  gives the following identity.
\begin{equation}\sum_{k=1}^{\infty}
\left ( \LL  \left ( \frac{ 1}{ \phi^2 k^2}\right )  +
\LL  \left ( \frac{1}{ \phi^4 k^2}\right )
+2
\LL  \left ( \frac{ 1}{ (\phi k -1)(\phi^2k-1)}\right )
\right ) = 2\pi^2/3.\end{equation}

\noindent  Restricting to rational vertices and $n=4$, one has $v_1=0$, $v_2= p/q$, $v_3= 1$, $v_4
= \infty$. Set $r= q-p$.
Then (13.2) reduces to the identity
\begin{equation}
\sum_{k=1}^{\infty}
\left ( \LL  \left ( \frac{ p^2}{ q^2 k^2}\right )  +
\LL  \left ( \frac{ r^2}{ q^2 k^2}\right )
+2
\LL  \left ( \frac{ pr}{ (qk-p)(qk-r)}\right )
\right ) = 2\pi^2/3.\end{equation}

More generally, if $n\ge 4$ and successive vertices are rational and Farey neighbors, the dilogarithm terms in (13.2) are all in terms of rationals with numerator $1$.

 \section{Continued fractions and cross ratios }
 The main purpose of this section is to give a technical lemma that will enable us to express the dilogarithm identities associated with continued fractions $\alpha$ of even period greater than two in terms of the convergents of the $k$th cyclic permutations $\alpha^{(k)}$ of $\alpha$. 

 \subsection{A technical lemma}
  Let $\alpha =[a_0,a_1,\cdots]$ be an infinite continued fraction and  $r_n= p_n/q_n$ its $n$th  convergent  as in subsection 9.1.
 Denote by $[r_i :r_j]$ the determinant of the matrix ${\tiny \left (
\begin{array}{cc}
p_i & p_{j} \\
q_i & q_{j}\\
\end{array}
\right ) }$. Note that 
\begin{equation}[r_i :r_j]=-[r_j :r_i].\end{equation}
Let $k$ be a fixed integer.
 Since the determinant function is a multi-linear function in terms of its
  columns, $d_k(n) :=[r_k : r_n]$
is a recurrence that shares the same recurrence with the $p_n$'s and $q_n$'s (see subsection 9.1).
 To be more precise, one has
\begin{equation}
d_k(k)  =0,\,\, d_k(k+1)= (-1)^{k-1},\,\,
d_k(n) =
a_{n}d_k(n-1)+d_k(n-2).\end{equation}


\begin{example} Let $\alpha = [\, \overline {1,2,3}\,]$ be a continued fraction
 of period 3. The
 first few convergents $r_k$'s are given as follows
$$
 r_0= \frac{1}{1},\,\,
 \frac{3}{2},\,\,
  \frac{10}{7},\,\,
  \frac{13}{9},\,\,
 \frac{36}{25}\,\,
  \frac{121}{84},\,\,
 \frac{157}{109},\,\,
 \frac{435}{302},\,\,
  \cdots
 $$
 Note that $a_0=1, a_1=2, a_2=3$.
In the case $k=0$, by (14.2), the first few $d_0(n)$'s  are
$$d_0(0)=0,\,\,d_0(1)=-1,\,\,-3\,\,-4,\,\,-11,\,\,-37,\,\,-48,\cdots.$$
In the case $k=2$,  one has $d_2(2+k)=-q_{k-1}$.
The first few $d_2(2+k)$'s are
$$d_2(2)=0,\,\,d_2(3)=-1,\,\,-2,\,\,-7\,\,-9,\,\,-25,\,\,-84,\,\,-109,\cdots.$$

\end{example}

\begin{lemma} Let $\alpha$, $a_i$ and  $d_n(m)$ be
 given as in $(9.1), \,(14.2) $ and let $r_n, r_{n+2}, r_{m}, r_{m+2}$
 be   convergents of $\alpha$, where $m \neq n$. Then

 \begin{equation}
 [ r_{n+2}, r_n, r_{m+2}, r_m]=
   \frac{ (-1)^{m+n}  a_{n+2}a_{m+2}}
   {d_m(n)d_{m+2}(n+2)}
   .\end{equation}
   
\end{lemma}
Note that $d_n(m)d_{n+2}(m+2)= d_m(n) d_{m+2}(n+2)$ by (14.1).
\begin{proof} Apply a  direct calculation and use (14.2) and the definition of $d_m(n)$.
	\end{proof}

More generally, one can determine $[r_a, r_b, r_c, r_e]$ as well but as we will see in equation (15.5), the restricted form in (14.3) suffices for our purpose.

 \section{Dilogarithm identities for continued fractions of even period $> 2$.}

Let $\alpha=[\overline{a_0,a_1,\cdots,a_{l-1}}]$ be a continued fraction of period $l$. Doubling the period
  if necessary, for example, 
   $\alpha=[\overline {1,2,3}]=[\overline {1,2,3,1,2,3}]$,
    we may assume that $l$ is even.
   Let
    $r_n= p_n/q_n$ be the $n$-th convergent of $\alpha$, and following standard convention, set $r_{-2}=0=0/1$, $r_{-1}=\infty=1/0$.
     One sees easily that
  \begin{equation} r_0 <r_2<r_4 <
 \cdots < \alpha <   \cdots < r_{l-1} < r_{l-3}<\cdots < r_3 <r_1.\end{equation}

\noindent
Since $\alpha$ has period $l$, letting  $A =
\left ( \begin{array}{cc}
p_{l-1}&  p_{l-2}\\
q_{l-1} &q_{l-2}\\
\end{array}
\right )$, a direct calculation shows that
   \begin{equation}
    A^n(r_k) =
\left ( \begin{array}{cc}
 p_{l-1}&  p_{l-2}\\
 q_{l-1} &q_{l-2}\\
\end{array}
\right )^n r_k = r_{nl+k}\,,\,\mbox{  where   } 0\le k\le l-1,\quad n\ge 0.\end{equation}

 \subsection{A feasible pair $(A, P)$ for $\alpha$} Let $\overline \alpha $
 the
  Galois conjugate of $\alpha$.
  Let $v_0=[\overline \alpha, \alpha]\cap[0, \infty]$ and 
  $v_1=[\overline \alpha, \alpha]\cap[r_{l-2}, r_{l-1}]$.
  Let  $P$ be the hyperbolic convex hull of
    \begin{equation}
    \{\infty, v_0,v_1, r_{l-1}, r_{l-3},\cdots, r_3, r_1\}.\end{equation}

    The sides $[v_0, \infty]$ and $[v_1, r_{l-1}]$ of $P$ are
 paired by $ A= {\tiny\left ( \begin{array}{cc}
 p_{l-1}&  p_{l-2}\\
 q_{l-1} &q_{l-2}\\
\end{array}
\right )}$, giving a crown $S$ with $l/2$ boundary cusps. Note that $A$ has determinant one  ($l$ is even) and that
 $(A, P)$ is a feasible pair.

\subsection{The sides of the universal cover $\tilde S$ of $S$}   Applying $A$ to the sides of $\tilde S$, the universal cover of $S$, and using (15.2),  we see that the sides of  $\tilde S$ split into $l/2  +1$
      orbits under the action of $A$. The side
      $[\overline \alpha, \alpha]$ is the only representative
        in its orbit. The other orbits are infinite sets.  The representatives of these orbits are
           \begin{equation}
   \{[\overline \alpha, \alpha]\}\cup   \{[r_1, r_{-1}]\} \cup  \{(r_3, r_1)\}
  \cup  \{[r_5, r_3]\} \cup \cdots \cup   \{[r_{l-1}, r_{l-3}]\} .
  \end{equation}

  \subsection{The dilogarithm identity for $\alpha$}
  Recall that $\bar \alpha$ and $\alpha$ are the fixed points of $A$. Similar to (13.2), and using the enlarged set of orthogeodesics $\hat O(S)$ of $S$ and applying Bridgeman's identity (13.1), one has

\begin{proposition}
  Let $\alpha=[\overline{a_0, a_1,a_1,\cdots,a_{l-1}}]$ be a continued
   fraction of even period $l$ and let $e_{2i+1}=[r_{2i+1}, r_{2i-1}]$, where $r_i=p_i/q_i$ is the $i$th convergent of $\alpha$. Then
   
   \begin{multline*}\mathop{\sum\sum}_{0 \le j<i\le l/2-1} \LL ([e_{2i+1}, e_{2j+1}])+
    \sum  _{j=0}^{l/2-1}
   \sum  _{i=0}^{l/2-1} \sum _{k=1}^{\infty}
     \mathcal L ([A^k(e_{2i+1}), e_{2j+1}])\\
   +
   \sum_{m=0}^{l/2-1}
   \mathcal L([\overline \alpha, \alpha, r_{2m+1}, r_{2m-1}])
   = \frac {\pi^2 l}{6}.
   \end{multline*}
    
    \end{proposition}

  We call the identity in Proposition 15.1  the dilogarithm identity associated to the continued fraction $\alpha$. The cross ratios occurring in the double and triple summations can be expressed in terms of the convergents of the $k$th cyclic permutations of $\alpha$ (see section 15.4) by applying Lemma 14.2. To start with, a direct application of Lemma 14.2 gives:

\begin{equation}
[A^k(e_{2i+1}), e_{2m+1}])
=
[r_{2n+1}, r_{2n-1}, r_{2m+1}, r_{2m-1}]
 = \frac{
 a_{2n+1}a_{2m+1}}
 {d_{2m-1}(2n-1)
 d_{2m+1}(2n+1)},\end{equation}

\smallskip
\noindent where $n=kl/2+i$,  since $A^k(e_{2i+1})= [r_{kl+2i+1}, r_{kl+2i-1}]$ by (15.2).
 Note that the numerator of (15.5) is bounded, 
  since $a_{l+k}=a_k$. In fact, if $a_{2k+1}=1$ for all $k$, then the numerator is always equal to 1. The denominator of (15.5) is a product of two
  recurrences that can be calculated easily from (14.2), see also Example 14.1. It can can also be expressed in terms of the convergents of the cyclic permutations of $\alpha$ as shown in the next subsection.

\subsection{The denominator of the cross ratio terms of Proposition 15.1}

Recall that $\alpha =[\overline {a_0. a_1, a_2, \cdots, a_{l-1}}]$ and $l$ is even. Define the $k$th (cyclic) permutation of $\alpha$ by
 \begin{equation}
 \alpha ^{(k)}:= [\overline {a_k, a_{k+1},\cdots,
  a_{l-1}, a_0, a_1, \cdots, a_{k-1}}].\end{equation}

  \noindent
   It is clear that $\alpha=\alpha^{(0)}$.
   Denote  the $s$-th convergent of $\alpha^{(k)}$ by $$p^{(k)}(s)/
   q^{(k)}(s).$$
   Note that we are writing $s$ as an argument instead of an index unlike the case of the convergents for $\alpha$, to fit with later applications, so $p^{(0)}(s)/q^{(0)}(s)=p_s/q_s$. 
   
   Suppose that $k$ is odd. By (14.2), one has
      \begin{equation}
      d_k(k)=0,\,\, d_k(k+1)=1,\,\, d_{k}(k+ t+2) = p^{(k+2)}(t).\end{equation}

The numerator of (15.5) is coming from the terms that define the continued fraction of $\alpha$. By (15.5) and (15.7), the denominator of the terms in the double and triple summation of Proposition 15.1
 can be described as follows.

\medskip
\begin{proposition} With terms as defined in Lemma 14.2, Proposition 15.1, and equation (15.5), we have:

\begin{equation}
\mathcal L \left (\frac{
 a_{2n+1}a_{2m+1}}
 {d_{2m-1}(2n-1)
 d_{2m+1}(2n+1)}\right )
 =\mathcal L \left (\frac{
a_{2n+1}a_{2m+1}}
 {p^{(2m+1)} (s)
 p^{(2m+3)}(s)}\right ),\end{equation}

\smallskip
\noindent where
$s = 2n-2m-2$ and
$p^{(k)}(s) $ is the numerator of the $s$-th convergent
 of the
  continued fraction of $\alpha^{(k)}$, the $k$th permutation of $\alpha$.
   
\end{proposition}

\smallskip

Combining Proposition $15.1$, $(15.5)$ and
$(15.8)$ gives the following dilogarithm identity
associated with $\alpha$:

\begin{theorem} Let $\alpha=[\overline{a_0,a_1, \ldots, a_{l-1}}]$ be a periodic continued fraction of even period $l$, $\bar{\alpha}$ its Galois conjugate, $\alpha^{(k)}=[\overline{a_k,a_{k+1},\ldots,a_{k-1}}]$ the $k$th cyclic permutation of $\alpha$, and $r^{k}(s)=p^{(k)}(s)/q^{(k)}(s)$ the $s$th convergent of $\alpha^{(k)}$. Then 
   
\begin{multline}\mathop{\sum\sum}_{0 \le j<i\le l/2-1} \LL \left(\frac{a_{2i+1}a_{2j+1}}{p^{(2j+1)}(2i-2j-2)p^{(2j+3)}(2i-2j-2)}\right)\\
+
\sum  _{j=0}^{l/2-1}
\sum  _{i=0}^{l/2-1} \sum _{k=1}^{\infty}
\mathcal L \left(\frac{a_{2i+1}a_{2j+1}}{p^{(2j+1)}(kl+2i-2j-2)p^{(2j+3)}(kl+2i-2j-2)}\right)\\
+
\sum_{m=0}^{l/2-1}
\mathcal L([\overline \alpha, \alpha, r_{2m+1}, r_{2m-1}])
= \frac {\pi^2 l}{6}. 
\end{multline}

\end{theorem}

\smallskip

\medskip

\section{Final remarks.} 
\begin{enumerate}
	\item 
Theorem 15.3 indicates  that
 if $\alpha $ has even period greater than 2, then the connection
  between  the arguments of
  $\mathcal L$ in the identity associated with $\alpha$  is not as tight as in (9.2). They can nonetheless be expressed in terms of the convergents of the $k$th permutations $\alpha^{(k)}$ of $\alpha$. Note that when $l=2$, the first double sum disappears, the triple sum reduces to a single sum involving only convergents of $\alpha$ and there is only one term in the finite sum. The identity is then equivalent to (9.2).

  \item  There are a couple of analogous identities associated to $\alpha$ which we can derive. Firstly, we can also consider the double crown associated to $\alpha$ (which has $l/2$ boundary cusps on both boundaries, so a total of $l$ boundary cusps) and obtain an identity from it analogous to (4.6), and which does not involve $\alpha$ explicitly. We leave details to the interested reader. Secondly, we may consider the other crown associated to $\alpha$ whose tines are the convergents of $\alpha$ with even indices.  Details are again left to the reader.

\item Finally, we remark that we have not considered here the case of finite identities derived from Bridgeman's theorem. As can be seen in \cite{B1} and \cite{B2}, many well-known identities for the Rogers dilogarithm follow easily from Bridgeman's identity (3.2) applied to various finite sided ideal hyperbolic polyhedra. We can for example consider the finite polygons in Section 13 with rational vertices, without a group action. Indeed, the literature on finite identities for the Rogers dilogarithm is quite extensive, see for example \cite{Kir}. The most interesting question is if it is possible to find exact values for $\LL(x)$ for exact arguments $x$, which are different from the small number of known cases: $\LL(0)$, $\LL(1/2)$, $\LL(1)$, $\LL(\phi^{-1})$ and $\LL(\phi^{-2}$), using some of these relations. There appear to be some evidence that there are theoretical algebraic obstructions (see \cite{Zag} for example) to finding new examples, so any new examples would be surprising.

\end{enumerate}






\medskip

\end{document}